\documentclass[final,times]{elsarticle}

\usepackage{latexsym,amssymb,amsmath,amsthm}

\usepackage{colortbl}

\usepackage[T1]{fontenc}
\usepackage{amsfonts}

\usepackage{verbatim}
\usepackage{latexsym,amssymb,amsmath}
\usepackage{amsfonts}
\usepackage{url}
\usepackage{hyperref}

\journal{Elsevier}

\begin{document}

\begin{frontmatter}

\title{Impact of  Delay on  Predator-Prey Models}

\author[adbel]{A. Moujahid}
\ead{abdelmalik.moujahid@ehu.eus}

\author[fernando]{F. Vadillo\corref{cor1}}
\ead{fernando.vadillo@ehu.eus}

\address{Department of Mathematics, University of the Basque Country}

\cortext[cor1]{Corresponding authors}

\begin{abstract}
Mathematical modeling based on time-delay differential equations is an important tool to study the role of delay in biological systems and to evaluate its impact on the asymptotic behavior of their dynamics. Delays are indeed found in many biological, physical, and engineering systems and are a consequence of the limited speed at which physiological, chemical, or biological processes are transmitted from one place to another. Since real biological systems are always subject to perturbations that are not fully understood or cannot be explicitly modeled, stochastic delay differential systems (SDDEs) provide a more realistic approximation to these models.
In this work, we study the predator-prey system considering three time delay models: one deterministic and two types of stochastic models. Our numerical results show relevant differences in their respective asymptotic behaviors.
\end{abstract}

 \bigskip

\end{frontmatter}

\bigskip

 \bigskip


\section{Introduction}\label{intro}
In recent years, life sciences have been the focus of numerous research activities to achieve a better understanding of many phenomena of life and the interaction between living systems and the environment.
In this context, mathematical modeling based on differential equations, especially stochastic differential equations with time delay, is an important tool for the study of biological and ecological systems. Mathematical models in population dynamics is a topic that began many years ago, and there is an extensive bibliography. Two classic references include. \cite[Chap. 1]{1973_May} from 1973 and also \cite[Chap. 1]{murray1} from 1989, this topic even appears in introductory courses on differential equations, (see, e.g., \cite[pp. 284]{simmons} or \cite[pp. 74]{boyce_diprima}). Four textbooks of interest are \cite{2006_devries}, \cite{2004_neal}, \cite{2014_matlab}, and \cite{2016_chou}, which help students explore and discover mathematical concepts and apply these concepts in building and analyzing mathematical models in life science disciplines such as biology, ecology, and environmental science.

However, many of these models ignore some hidden processes and external influences that are poorly understood and generally affect the dynamics of these systems. Therefore, for a more realistic interpretation and approximation to real systems, it is necessary to consider stochastic delay differential equations \cite{2003_allen,2007_allen,2009_chasnow,2016_pardoux} and \cite{2022_rihan}, where the relevant mechanisms are modeled as stochastic processes.

On the other hand, delayed interacting processes are ubiquitous in many biological systems, which generally consist of interacting units, and depending on the time scales of these units, the speed of information propagation becomes relevant to their dynamics. Therefore, delayed stochastic differential equations are crucial for the study and better understanding of these systems. Moreover, it is almost self-evident that the introduction of memory terms greatly increases the complexity.

The rest of the paper is organized as follows. In Section \ref{s:deter}, we present the deterministic model for the prey-predator system and study its asymptotic behavior as a function of time delay. Section \ref{s:stochastic1} describes the first delayed stochastic model constructed by adding Gaussian noise to the deterministic terms. The second delayed stochastic model, based on probabilistic considerations, is derived in Section \ref{s:stochastic2}.
we conclude in Section \ref{conclusiones} and draw the main contributions as well as possible future applications to more complex systems with more parameters.


\bigskip
\section{The Prey-Predator Deterministic Model}\label{s:deter}

The predator-prey model was originally proposed by A. J. Lotka and V. Volterra in the 1920s. Since then, many other models have been proposed with different functional responses. In this paper, we consider the Michaelis-Menten type predator-prey model introduced by Freedman \cite{Freedman1980} in 1980. The time-delay dynamics of this model is governed by the following system of differential equations,
\begin{equation}\label{e:dde}
\left\{ \begin{array}{l}
    \dfrac{dx}{dt} ~=~    r~x(t) \left( 1 - \dfrac{x(t)}{K } \right)   ~- ~          \dfrac{\beta ~ x(t) ~ y (t-\tau) }{1+ \sigma ~ x(t)}, \\    \\[0.5mm]
    \dfrac{dy}{dt} ~=~ \dfrac{\beta ~ x(t) ~ y (t-\tau) }{1+ \sigma ~ x(t)} ~- ~ a ~y(t), \\
     \end{array} \right.
\end{equation}
where $x(t)$ and $y(t)$ represent the population densities of the prey and predator, respectively. The prey population is assumed to have logistic growth with a carrying capacity $K$ and a specific growth rate constant $r$. The feeding rate $\beta$ is the maximum number of prey that can be eaten by a predator in any unit of time, $r$ is a specific growth rate, $\sigma$ is a positive constant describing the effects of capture rate, and $a$ refers to the dead rate for the predator.

The initial conditions are set to $x(0)>0$ and $ y (t)=\psi(t)>0$ when $t \in [ -\tau, 0]$, where $\psi(t)$ is a smooth function and $\tau$ the time delay. In \cite[pp. 217]{2021_rihan} it is shown that there are three fixed points: the trivial equilibrium $\varepsilon_0=(0,0)$, the semi-trivial equilibrium $\varepsilon_1=(K,0)$, and the interior equilibrium $\varepsilon_+ ~=~ (x^*,y^*)$, with
  \begin{equation*}
 x^* ~=~ \dfrac{a}{\beta- \sigma a}, \quad  y^* ~=~  \dfrac{r (x^*)^2}{Ka} (R_0-1),  \quad \hbox{and} \quad \mathcal{R}_0 ~=~  \dfrac{K}{x^*}.
\end{equation*}

\noindent The stability theorem \cite[Th. 11.2]{2021_rihan} of this interior equilibrium shows that $\varepsilon_+ ~=~ (x^*,y^*)$ is locally stable when $ 1 < \mathcal{R}_0 \le \mathcal{R}_c=~ 2+  \dfrac{1}{1+2 \sigma x^*} $, while for $ 1 < \mathcal{R}_c < \mathcal{R}_0 $, there exist a bifurcation $\tau^*>0$ such that the equilibrium is asymptotically stable for values of $\tau \in [0, \tau^*]$, and unstable when $\tau > \tau^*$.

\begin{figure}[hbt!]
    \centering
    \includegraphics[width=.8\textwidth]{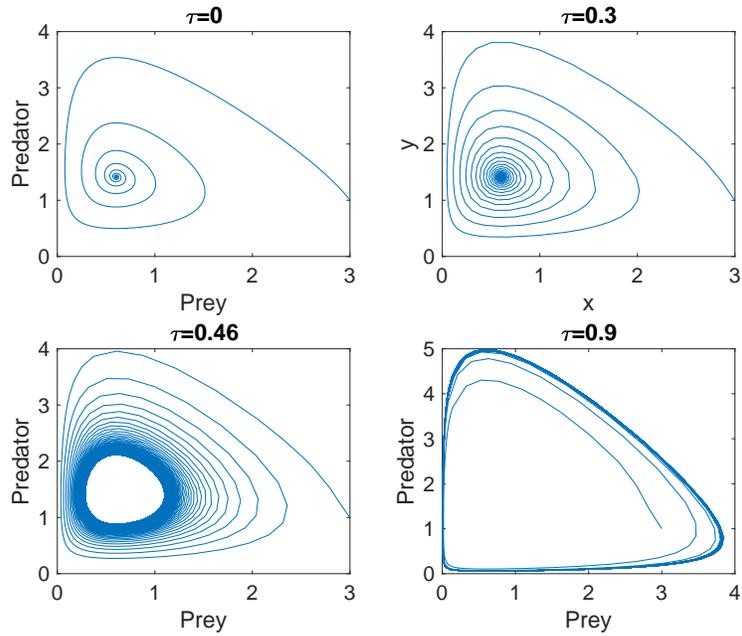}
    \caption{ Phase portraits of the model \ref{e:dde} at different time-delays. $r = 0.8, K = 5, \sigma = 0.01, \beta = 0.5, a = 0.3 \Rightarrow \tau^*=0.46$ and $0 \le t \le 10000$. $x(0)=3, y(0)=1$}.
    \label{fig1}
\end{figure}

Figure \ref{fig1} shows the phase portraits of the deterministic system given by Eq. \ref{e:dde} for different values of the time-delay $\tau$, and for the parameter values  $r = 0.8, \quad  K = 5, \quad  \sigma = 0.01, \quad  \beta = 0.5, \quad  a = 0.3
$ \cite{shampine}. For these values, the interior equilibrium is $(x^*,y^*)=(0.6036,1.4153)$, $\mathcal{R}_0= 8.2833$ and $\mathcal{R}_c=2.9881$. This corresponds to the case where $1 < \mathcal{R}_c< \mathcal{R}_0$. As can be seen, there is a critical value of $\tau$, the system switches from a stable fixed point to a limit cycle.

\begin{figure}[hbt!]
    \centering
    \includegraphics[width=\textwidth, height=5cm]{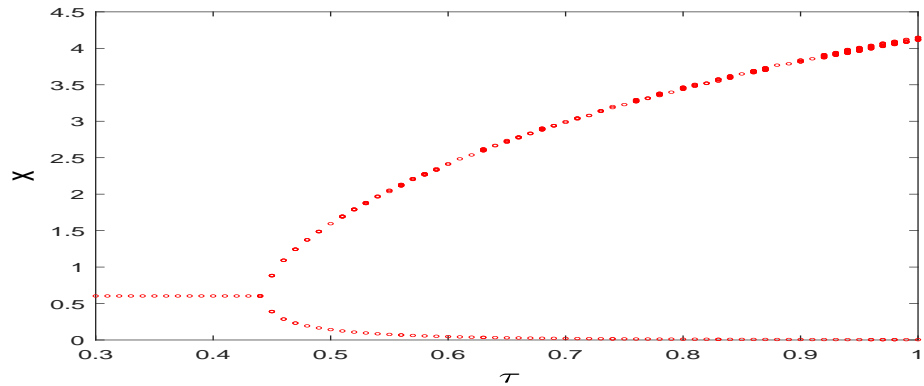}
    \caption{ Hopf bifurcation diagram of the deterministic system \ref{e:dde} considering the time delay as a bifurcation parameter. For $\tau<0.46$, the interior equilibrium is asymptotically stable. For $\tau\geq0.46$ an asymptotically limit cycle emerges.}
    \label{fig2}
\end{figure}

To better visualize these bifurcation patterns, we also plotted the bifurcation diagram with $\tau$ as a varying parameter (see Figure \ref{fig2}). According to our numerical results, the interior equilibrium is asymptotically stable for values of $\tau$ smaller than the bifurcation value $\tau^+=0.46$ at which the Hopf bifurcation occurs. At this critical value, the trajectories of the system slowly converge toward the limit cycle, while for values of $\tau$ larger than $\tau^+$, the behavior of the system quickly collapses into the limit cycle (see Figure \ref{fig1}).

\bigskip

\section{White noise based stochastic model: Model 1}\label{s:stochastic1}


From a biological perspective, the inclusion of stochastic perturbations in the deterministic predator-prey model \eqref{e:dde} allows for a more realistic interpretation of the system. The simplest way to do this is to add a white noise that is proportional to $x(t)$ and $y(t)$, leading to the SDDE:

\begin{equation}\label{model+1}
\left\{ \begin{array}{l}
             d x(t) = \left[  r~ x(t) \left( 1 - \dfrac{x(t)}{K } \right)- \dfrac{\beta ~ x(t) ~ y(t-\tau)}{1+ \sigma ~ x(t)} \right] ~ dt +  \nu_1 x(t) ~d W_1(t),\\[4mm]
              d y(t) = \left[ \dfrac{\beta ~ x(t) ~ y(t-\tau)}{1+ \sigma ~ x(t)} - a y(t) \right] ~ dt +  \nu_2 y(t) ~d W_2(t),
           \end{array} \right.
\end{equation}
where $W_1(t)$ and $W_2(t)$  denote two  independent Brownian motions, and the positive constants $\nu_1$ and $\nu_2$ refer to the intensities of the white noises.

The numerical simulations use Milstein's scheme for SDDEs, discussed in \cite{2006_milstein} or \cite[p. 137]{2021_rihan}, actually it is a simple adaptation of the method for SDE of the same name
 which can be found in classical texts, e.g. \cite{2001_Higham} , \cite{kloeden}, and more recently \cite{2021_high_kloeden}.

\begin{figure}[hbt!]
    \centering
    \includegraphics[width=\textwidth]{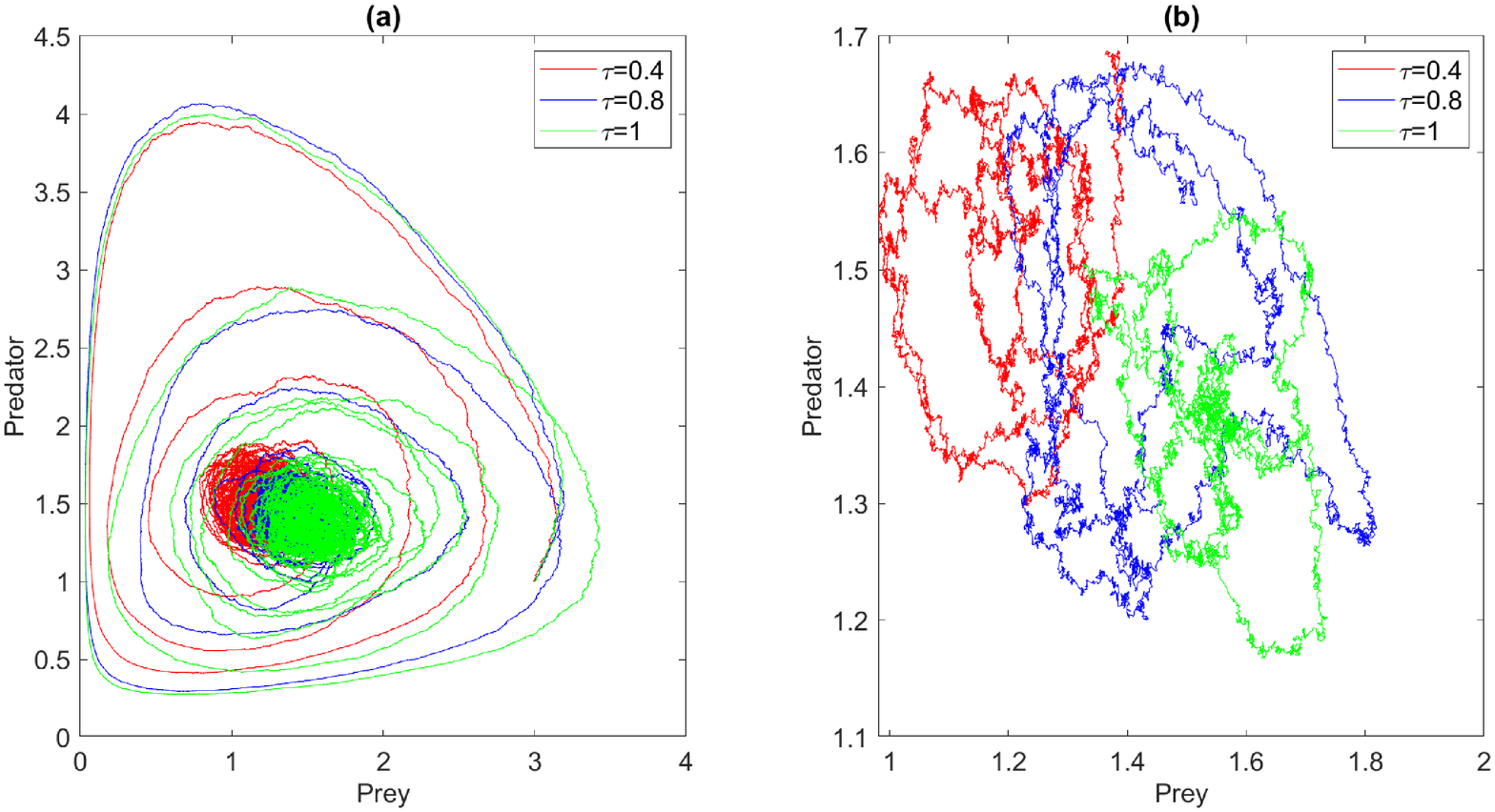}
    \caption{Model 1: Phase portraits of the model for different values of the time delay $\tau$ for $\nu_1=\nu_2=0.1$, $r=0.8$, $K=5$, $\sigma=1/3$, $\beta=0.5$ and $a=0.3$. The results were obtained from the stochastic data over 500 simulations and $0 \le t \le 1500$. For the phase portrait in panel (b), we considered only the permanent regime.}
    \label{fig3}
\end{figure}

In Figure \ref{fig3} we plotted the phase portraits the model for different values of the time delay $\tau$. The results were obtained from the stochastic data over 500 simulations and $0 \le t \le 1500$. As you can see (panel (a)), for each value of $\tau$, the trajectories approach different fixed points, and no Hopf bifurcation is observed. In other words, the time delay move the system to different regions of the state space where a particular stationary fixed point occurs.

\begin{figure}[h!]
    \centering
    \includegraphics[width=0.9\textwidth, height=5cm]{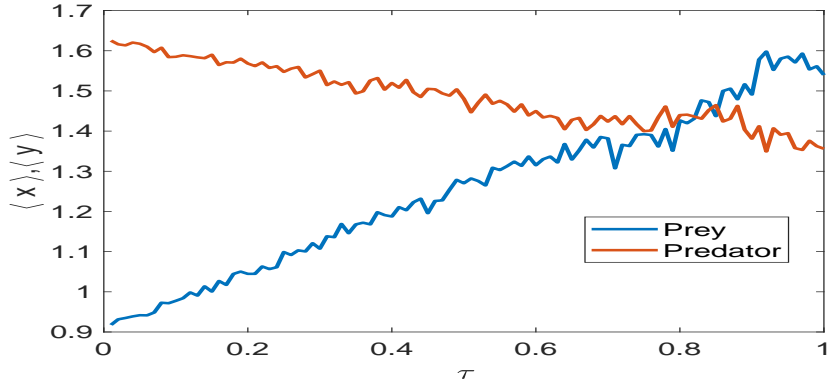}
    \caption{Model 1: The long-term average evolution of the stationary equilibrium as a function of the time delay. $\nu_1=\nu_2=0.1$, $r=0.8$, $K=5$, $\sigma=1/3$, $\beta=0.5$ and $a=0.3$.}
    \label{fig4}
\end{figure}

\begin{figure}[hbt!]
    \centering
    \includegraphics[width=\textwidth]{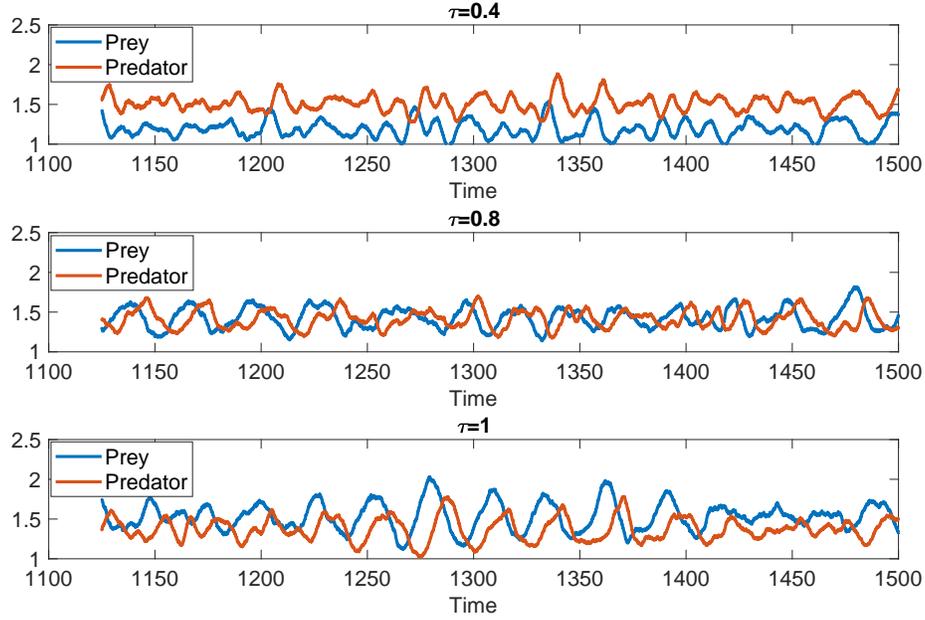}
    \caption{Model 1: Long-term trajectories of the model for different values of the time delay $\tau$ for $\nu_1=\nu_2=0.1$, $\nu_1=\nu_2=0.1$, $r=0.8$, $K=5$, $\sigma=1/3$, $\beta=0.5$ and $a=0.3$. The results were obtained from the stochastic data over 500 simulations and $0 \le t \le 1500$.}
\label{fig5}
\end{figure}

To better understand how the time delay affects the system dynamics, we also plotted the long-term average evolution of the steady-state equilibrium as a function of the time delay (see Figure \ref{fig4}). For each value of $\tau$, the trajectories were averaged over 500 runs and $0 \le t \le 1500$.
For values of $\tau$ less than $\tau{**}=0.8$, the prey population is on average smaller than the predator population, while for values greater than $\tau{**}$ the situation reverses and the predator population dominates the prey population. This value of $\tau^{**}$ is of course dependent on the capture rate $\sigma$, and increases with decreasing $\sigma$.
At the critical value $\tau{**}$, the population sizes of prey and predators are equal on average (see Figure \ref{fig5}).

\section{Probabilistic stochastic model: Model 2}\label{s:stochastic2}


The second stochastic model that we will consider, which we could call the probabilistic model, is very similar to the ones analyzed
in \cite{2012_hoz_vadillo},  \cite{2019_lotka}, \cite{sis-model} and \cite{2021_abdel_me}. Assuming that the changes and their probabilities are those given in Table 1.

\begin{table}[!ht]
\begin{center}
\caption{Possible change in $\mathbf{Z}=(X , Y)^T$ and their probabilities }
\begin{tabular}{l|l} \label{t:1}
   Change  &  Propability  \\
   \hline
         &  \\
         $\Delta \mathbf{Z}^{(1)}~=~ (1, 0  )^T $ &  $ p_1=  r~ x(t) \left( 1 - \dfrac{x(t)}{K } \right) ~ \Delta t $ \\
         &  \\
         $\Delta \mathbf{Z}^{(2)}~ =~(-1 , 0)^T $ &  $ p_2=  \dfrac{\beta ~ x(t) ~ Y(t-\tau)}{1+ \sigma ~ x(t)} ~ \Delta t $ \\
         & \\
         $\Delta \mathbf{Z}^{(3)}~=~ (0, 1)^T $ &  $ p_3 =  \dfrac{\beta ~ x(t) ~ Y(t-\tau)}{1+ \sigma ~ x(t)} ~ \triangle t $ \\
         & \\
         $\Delta \mathbf{Z}^{(4)}~ =~ ( 0, -1 )^T $ & $ p_4=  a y(t)   ~ \Delta t $ \\
\end{tabular}
\end{center}
\end{table}

\noindent Fixing $\mathbf{Z}(t)$ at time $t$, we calculate the expected change for $\mathbf{Z}=(X, Y)^T$
\begin{equation}
   \mathrm{E}( \Delta \mathbf{Z}) ~=~  \sum_{j=1}^{4} p_j  ~ \Delta \mathbf{Z}^{(j)} ~=~
    \left(\begin{array}{c}
        r~ x(t) \left( 1 - \dfrac{x(t)}{K } \right)- \dfrac{\beta ~ x(t) ~ y(t-\tau)}{1+ \sigma ~ x(t)}, \\
        \dfrac{\beta ~ x(t) ~ y(t-\tau)}{1+ \sigma ~ x(t)} - a y(t),
   \end{array} \right) ~\Delta t,
\end{equation}

\noindent and the covariance matrix
\begin{equation}
   \mathrm{E}( \Delta \mathbf{Z} ( \Delta \mathbf{Z})^T) ~=~  \sum_{j=1}^{4} p_j ~ (\Delta \mathbf{Z}^{(j)})  (\Delta \mathbf{Z}^{(j)})^T ~=~ \mathbf{D}(x,y)
     ~ \Delta t,
\end{equation}
\noindent where $\textbf{D}(x,y)$ is the diffusion matrix
\begin{align*}
 \mathbf{D} ~=~ \left(
                            \begin{array}{cc}
                              r~ x(t) \left( 1 - \dfrac{x(t)}{K } \right)+ \dfrac{\beta ~ x(t) ~ y(t-\tau)}{1+ \sigma ~ x(t)} & 0 \\
                             0 &  \dfrac{\beta ~ x(t) ~ y(t-\tau)}{1+ \sigma ~ x(t)} + a y(t)
                            \end{array}
                          \right),
\end{align*}
\noindent so we arrived at stochastic differential system
\begin{equation}\label{model+2}
\left\{ \begin{array}{l}
             d x(t) = \left[  r~ x(t) \left( 1 - \dfrac{x(t)}{K } \right)- \dfrac{\beta ~ x(t) ~ Y(t-\tau)}{1+ \sigma ~ x(t)} \right] ~ dt
             +  \sqrt{ r~ x(t) \left( 1 - \dfrac{x(t)}{K } \right)+ \dfrac{\beta ~ x(t) ~ y(t-\tau)}{1+ \sigma ~ x(t)}} ~d W_1(t),\\[4mm]
              d y(t) = \left[ \dfrac{\beta ~ x(t) ~ y(t-\tau)}{1+ \sigma ~ x(t)} - a y(t) \right] ~ dt
              +  \sqrt{\dfrac{\beta ~ x(t) ~ y(t-\tau)}{1+ \sigma ~ x(t)} + a y(t)} ~d W_2(t)
           \end{array} \right.
\end{equation}

\begin{figure}[h!]
    \centering
    \includegraphics[width=\textwidth]{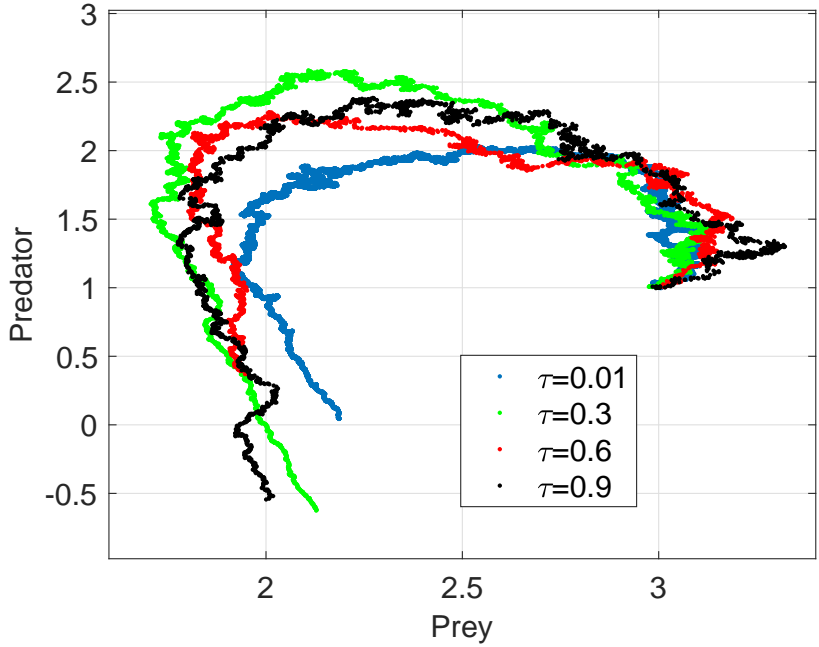}
    \caption{Model 2: mean of 100 simulations . For  $\tau=0.01$ in blue, $\tau=0.3$ green, $\tau=0.6$ in red with $0 \le t \le 9$ and $\tau=0.9$ in black with $0 \le t \le 12 $ } .
    \label{fig6}
\end{figure}

\begin{figure}[h!]
    \centering
    \includegraphics[width=\textwidth, height=4cm]{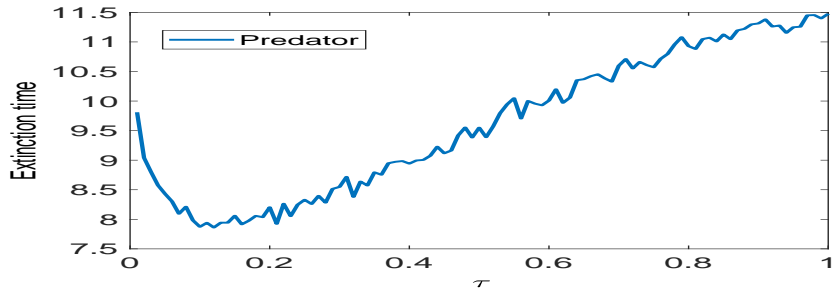}
    \caption{Model 2: Extinction time of the predator population as a function of time delay. The parameters of the model are: $r=0.8$, $K=5$, $\sigma=1/3$, $\beta=0.5$ and $a=0.3$ } .
    \label{fig7}
\end{figure}

In Figure \ref{fig6} we have plotted the phase portraits of Model 2 for different values of $\tau$, ranging from low to higher values. The results correspond to an average of 100 trials with initial value $x(0)=3, y(0)=1$ and $0 \le t \le 9$. In this case, the predator population goes to extinction very quickly, and the time to extinction appears to increase with delay, showing a minimum for low values of $\tau$ (see Figure \ref{fig7}).

\section{Conclusion}\label{conclusiones}
We have analysed the behaviour of two time-delayed stochastic differential systems for the prey-predator model and compared them with the corresponding deterministic system. The introduction of stochastic perturbations into these time-delayed models has a significant impact on population extinction and enriches the dynamics of the models. We found a significant difference between the deterministic system and the two classes of stochastic systems. Considering the time delay as a bifurcation parameter, the deterministic system exhibits a Hopf bifurcation leading to a critical $\tau^+$ from which the system transitions from the stable asymptotic fixed point to a stable limit cycle.

When the stochastic perturbation is modelled by simply adding a white noise to the deterministic components of the system (Model 1), the time delay appears to drive the system into different regions of state space where a particular stationary fixed point occurs, but no Hopf bifurcation is observed. Nevertheless, the population appears to be persistent for both the deterministic system and the stochastic model 1. However, in contrast to what we observed previously, for the second stochastic system (Model 2) derived from a probabilistic approach, the system collapses rapidly and the population goes to extinction. This is consistent with previous work supporting the fact that the nature of the stochastic system has a significant impact on the asymptotic behavior of the population, and environmental noise can cause the solution of the stochastic system to be extinct \cite {2019_lotka}.

\section*{Acknowledgments}
This work was supported by Spanish Ministry of Sciences Innovation and Universities with the project PGC2018-094522-B-100 and by the Basque Government with the project IT1247-19.

\bibliographystyle{amsplain}
\bibliography{ref_ddl_predator}

\end{document}